\documentclass[a4paper,12pt]{article}
\usepackage[T2A]{fontenc}
\usepackage[utf8]{inputenc} 
\usepackage{amssymb,amsfonts}
\usepackage{amsfonts,amsmath,amssymb,amscd,latexsym}
\usepackage{srcltx}

\linethickness{0.4pt}
\ifx\plotpoint\undefined\newsavebox{\plotpoint}\fi 

\newcounter{ppp}
\newcommand{\bi}{\bibitem}
\newcommand{\nb}{\newblock}

\newcommand{\be}[1]{\begin{equation}\label{#1}}
\newcommand{\ee}{\end{equation}}

\newcommand{\la}{\langle\,}
\newcommand{\ra}{\,\rangle}

\newcommand{\prf}{{\bf Proof.}\ }

\newcommand{\hgt}{\mathop{\rm ht}}

\newtheorem{thm}{\quad Theorem}
\newtheorem{lm}{\quad Lemma}

\newtheorem{df}{\quad Definition}
\newtheorem{prop}{\quad Proposition}

\title{Cayley graphs of R.\,Thompson's group $F$:\\ new estimates for the density}
\author{\vspace{2ex}
Victor Guba\thanks{This work is supported by the Russian Foundation
for Basic Research, project no. 20-01-00465.}\\
Vologda State University,\\
15 Lenin Street,\\
Vologda\\
Russia\\
160600\\
E-mail: gubavs{@}vogu35.ru}
\date{}

\begin{document}

\maketitle

\begin{abstract}

By the density of a finite graph we mean its average vertex degree. For the Cayley graph of a group $G$ with $m$ generators, it is known that $G$ is amenable if and only if the supremum of densities of its finite subgraphs has value $2m$.

For R. Thompson’s group $F$, the problem of its amenability is a long-standing open question. There were several attempts to solve it in both directions. For the Cayley graph of $F$ in the standard set of group generators $\{x_0,x_1\}$ there exists a construction due to Jim Belk and Ken Brown. It was presented in 2004. This is a family of finite subgraphs whose densities approach $3{.}5$. Many unsuccessful attempts to improve this estimate led to conjecture that this construction was optimal. This would imply non-amenabilty of $F$.

Recently we got an improvement showing that this conjecture turned out to be false. Namely, there exist finite subgraphs in the Cayley graph of $F$ in generators $x_0$, $x_1$ with density strictly exceeding $3{.}5$. This makes amenability of $F$ more truthful. Besides, we disprove one more conjecture showing that there are finite subgraphs in the Cayley graph of $F$ in generators $x_0$, $x_1$, $x_2$ with density strictly exceeding $5$.

\end{abstract}
\vspace{5ex}

Here is a brief description of the contents of our paper. In Section~\ref{prel} we recall some basic concepts we need. This includes Thompson's group $F$ and its positive monoid $M$; graphs in the sense of Serre and Cayley graphs of groups; density and isoperimetric constants; rooted binary trees and forests. 

In Section~\ref{sgpd} we describe explicitly the way of representing elements in $F$ by (non-sperical) semigroup diagrams together with the representaion of elements in $M$ by marked rooted binary forests. Also we explain how the generators of the group $F$ act on them.

In Section~\ref{dens} we discuss the amenability problem for the group $F$ in connection with its Cayley graph density. We list all known results in this direction including the one due to Belk and Brown. We formulate there one of our main results (Theorem~\ref{over3.5}) that the density of the Cayley graph of $F$ in standard generators $x_0$, $x_1$ stricltly exceeds $3{.}5$. This disproves some known conjectures.

In Section~\ref{gf} we introduce some of the most important technical tools for the proof of our results. Here we consider generating functions for various sets of marked rooted binary forests. 

Section~\ref{thm1} contains the proof of Theorem~\ref{over3.5}. In Section~\ref{thm2} we state and prove Theorem~\ref{dens012} showing that the density of the Cayley graph of $F$ in generators $x_0$, $x_1$, $x_2$ strictly exceeds $5$. This disproves one more of our conjectures from~\cite{Gu22}.

\section{Preliminaries}
\label{prel}

\subsection{Thompson's monoid and group}
\label{mf}

Let $X=\{x_0,x_1,x_2,\ldots,x_m,\ldots\}$ be an infinite countable alphabet. By $M$ we denote the monoid given by the following presentation

\be{xinf}
\la X\mid x_j{x_i}=x_ix_{j+1}\ (0\le i < j)\,\ra.
\ee

Applying rewriting rules of the form $x_jx_i\to x_ix_{j+1}$, where $j>i\ge0$, we get a word of the form $x_{i_1}\ldots x_{i_k}$, where $k\ge0$ and $0\le i_1\le\cdots\le i_k$. It is well-known that the result does not depend on the order of applying the rewriting rules. The corresponding word over $X$ is called the {\em normal form} of a monoid element.
\vspace{1ex}

It is easy to check that the monoid $M$ is cancellative. Also it is known that for any $a,b\in M$ there exists their least common right multiple. A classical Ore theorem states that for a cancellative monoid $M$ with common right multiples, there exists a natural embedding of $M$ into its group of quotients (see~\cite{Ljap} for details). Any element of this group belongs to $MM^{-1}$, and the group is given by the same presentation~(\ref{xinf}). It is denotet by $F$.
\vspace{1ex}

This group was introduced by Richard J. Thompson in the 60s. Details can be found in the survey \cite{CFP}; see also papers~\cite{BS,Bro,BG}. It is easy to see that $x_n=x_0^{-(n-1)}x_1x_0^{n-1}$ for any $n\ge2$, so the group is generated by $x_0$, $x_1$. This generating set will be called {\em standard}. The group $F$ can be given by the following presentation with two defining relations

\be{x0-1}
\la x_0,x_1\mid x_1^{x_0^2}=x_1^{x_0x_1},x_1^{x_0^3}=x_1^{x_0^2x_1}\ra
\ee
of length $10$ and $14$ respectively, where $a^b\leftrightharpoons b^{-1}ab$. 

Elements of the monoid $M$ are called {\em positive elements}. Each element of $F$ can be uniquely represented by a group {\em normal form\/}, that is, a word of the form $pq^{-1}$, where $p$, $q$ are normal forms of positive elements, 
and the following is true: if the word $pq^{-1}=\ldots x_i\ldots x_i^{-1}\ldots$ contains both $x_i$ and $x_i^{-1}$ for some $i\ge0$, then it also contains $x_{i+1}$ or $x_{i+1}^{-1}$.
\vspace{1ex}

Many authors use an equivalent definition of $F$ in terms of piecewise-linear functions. In fact, there are various presentations of $F$ in these terms. In our paper we are not going to use these functions at all so we do not list all conditions for them. This information is standard and can be found in may sources including~\cite{CFP,BS}. 
\vspace{1ex}

It is also known that $F$ is a {\em diagram group} over the simplest semigroup presentation ${\cal P}=\langle x\mid x=x^2\rangle$. See~\cite{GbS} for the theory of these groups. In our paper, we will use semigroup diagrams to represent elements of $F$, as well as marked binary forests for the same purposes. These techniques will be explained below.

\subsection{Graphs and automata}
\label{grr}

We will often work with graphs, including Cayley graphs of groups. Besides, we will use semigroup diagrams, (rooted binary) trees and forests. So we are going to recall the concept of a graph in the sense of Serre~\cite{Se80} together with the corresponding notation. 
Informally, this is non-oriented graph, where each geometric edge consists of two mutually inverse directed edges. 

Formally, this is a 5-tuple $\Gamma=\la V,E,^{-1},\iota,\tau\ra$, where $V$, $E$ are disjoint sets (of {\em vertices} and (directed) {\em edges}, respectively); $^{-1}$ denotes the mapping from $E$ to itself; $\iota$ and $\tau$ are mappings from $E$ to $V$. We assume that $e\ne e^{-1}$, $(e^{-1})^{-1}=e$, $\iota(e^{-1})=\tau(e)$, $\tau(e^{-1})=\iota(e)$ for all $e\in E$. Here $e^{-1}$ is called the {\em inverse} edge of $e$; $\iota(e)$ is the {\em initial} vertex of $e$; $\tau(e)$ is the {\em terminal} vertex of $e$.

A {\em path} in a graph $\Gamma$ is defined in a standard way. This is either a single vertex (called a {\em trivial} path), or a sequence of edges written as $p=e_1\ldots e_n$, where $\tau(e_i)=\iota(e_{i+1})$ for all $1\le i < n$.
\vspace{1ex}

Let $G$ be a group generated by a set $A$. Its {\em right Cayley graph} $\Gamma_r={\cal C}(G;A)$ is defined as follows. The set of vertices is $G$; for any $a\in A^{\pm1}$ we put a directed edge $e=(g,a)$ from $g$ to $ga$. The letter $a$ is called the {\em label} of this edge. The inverse directed edge $e^{-1}$ labelled by $a^{-1}$ goes from $ga$ to $g$. Labels are naturally extended to the set of paths in Cayley graphs.

Sometimes it is preferable to work with {\em left} Cayley graphs $\Gamma_l={\cal C}(G;A)$. They are defined in a similar way. The set of its vertices is also $G$. A directed egde $e=(a,g)$ with label $a$ goes from $ag$ to $g$. That is, going along the edge labelled by $a$, means cancelling $a$ on the left. Inverse edges and labels of paths are defined similarly for this case. Notice that a word $w$ in group generators will be the label of a path from the vertex $g\in G$ represented by $w$ to the identity element in case of left Cayley graphs.
\vspace{1ex}

The cardinality of a finite set $Y$ will be denoted by $|Y|$. 

Let $G$ be a group generated by $A$. For our needs we can assume that $A$ is always finite, $|A|=m$.  Let $\Gamma={\cal C}(G,A)$ be the Cayley graph of $G$, right or left. To any finite nonempty subset $Y\subset G$ we assign a subgraph in $\Gamma$ adding all edges connecting vertices of $Y$. So given a set $Y\subset G$, we will usually mean the corresponding subgraph. This is a labelled graph that we often call an {\em automaton}. For each $g\in Y$ we have exactly $2m$ directed edges in $\Gamma$ starting at $g$, where $a\in A^{\pm1}$. If the endpoint of such an edge with label $a$ belongs to $Y$, then we say that the vertex $g$ of our automaton $Y$ {\em accepts} $a$. For the case of right Cayley graphs this means $ga\in Y$, for the case of left Cayley graphs this means $a^{-1}g\in Y$. 

A vertex $g\in Y$ is called {\em internal} whenever it accepts all elements $a\in A^{\pm1}$. That is, the degree of $g$ in $Y$ equals $2m$. Otherwise we say that $g$ belongs to the {\em inner boundary} of $Y$ denoted by $\partial Y$. If a vertex does not belong to $Y$ but it is connected by an edge to a vertex in $Y$, then we say that the vertex belongs to the {\em outer boundary} of $Y$ denoted by $\partial_oY$.

An edge $e$ is called {\em internal} whenever it connects two vertices of $Y$. If a directed edge $e$ connects a vertex of $Y$ with a vertex outside $Y$, then we call it {\em external}. The set of external edges form the {\em Cheeger boundary} of $Y$ denoted by $\partial_{\ast}Y$. It is easy to check that the cardinalty of $\partial_{\ast}Y$ does not exceed cardinalities of $\partial Y$ and $\partial_oY$.
\vspace{1ex}

By the {\em density} of a finite subgraph $Y$ we mean its average vertex degree. This concept was introduced in~\cite{Gu04}; see also~\cite{Gu21a}. It is denoted by $\delta(Y)$. A {\em Cheeger isoperimetric constant} of the subgraph $Y$ is the quotient $\iota_*(Y)=\frac{|\partial_{\ast}Y|}{|Y|}$. It follows directly from the definitions that $\delta(Y)+\iota_*(Y)=2m$. Indeed, each vertex $v$ has degree $2m$ in the Cayley graph $\Gamma$. This is the sum of the number of internal edges starting at $v$, which is $\deg_Y(v)$, and the number of external edges starting at $v$. Taking the sum over all $v\in Y$, we have $2m|Y|$, which is equal to $\sum\limits_v\deg_Y(v)+|\partial_{\ast}Y|$. Dividing by $|Y|$, we get the above equality.
\vspace{1ex}

\subsection{Rooted binary trees and forests}
\label{rbtf}

We add this short subsection to introduce some notation used in the paper. More details on the subject can be found in~\cite{Gu21a}.

Formally, a {\em rooted binary tree} can be defined by induction.

1) A dot $\cdot$ is a rooted binary tree.

2) If $T_1$, $T_2$ are rooted binary trees, then $(T_1\hat{\ \ }T_2)$ is a rooted binary tree.

3) All rooted binary trees are constructed by the above rules.
\vspace{1ex}

Instead of formal expressions, we will use their geometric realizations. A dot will be regarded as a point. It is called a {\em trivial} tree. It coincides with its root. If $T=(T_1\hat{\ \ }T_2)$, then we draw a {\em caret\/} for $\,\hat{}\,$ as a union of two closed intervals $AB$ (goes left down) and $AC$ (goes right down). The point $A$ is the {\em root} of $T$. After that, we draw trees for $T_1$, $T_2$ and attach their roots to $B$, $C$ respectively in such a way that the trees have no intersection. 

\refstepcounter{ppp}
\begin{center}
\begin{picture}(115.833,52.721)(0,0)
\put(2,2){\line(1,1){18}}
\put(2,2.21){\line(1,1){18}}
\put(67.38,31.221){\line(1,1){18}}
\put(20,20){\line(1,-1){18.25}}
\put(20,20.21){\line(1,-1){18.25}}
\put(85.38,49.221){\line(1,-1){18.25}}
\put(19.75,23.5){\makebox(0,0)[cc]{$A$}}
\put(19.75,23.71){\makebox(0,0)[cc]{$A$}}
\put(85.13,52.721){\makebox(0,0)[cc]{$A$}}
\put(6.25,2.5){\makebox(0,0)[cc]{$B$}}
\put(6.25,2.71){\makebox(0,0)[cc]{$B$}}
\put(71.63,31.721){\makebox(0,0)[cc]{$B$}}
\put(33.426,2.102){\makebox(0,0)[cc]{$C$}}
\put(33.426,2.312){\makebox(0,0)[cc]{$C$}}
\put(98.805,31.323){\makebox(0,0)[cc]{$C$}}
\multiput(67.482,30.903)(-.0337322344,-.0566219649){349}{\line(0,-1){.0566219649}}
\multiput(55.709,11.142)(3.333554,.030032){7}{\line(1,0){3.333554}}
\multiput(79.044,11.352)(-.0337094044,.0588382331){343}{\line(0,1){.0588382331}}
\multiput(103.43,31.323)(-.0336722591,-.0418628086){462}{\line(0,-1){.0418628086}}
\put(87.874,11.983){\line(1,0){27.96}}
\multiput(115.833,11.983)(-.0337257921,.0522749777){374}{\line(0,1){.0522749777}}
\put(67.272,19.341){\makebox(0,0)[cc]{$T_1$}}
\put(102.589,19.341){\makebox(0,0)[cc]{$T_2$}}
\put(15.346,33.426){\makebox(0,0)[cc]{a caret}}
\put(79.465,3.784){\makebox(0,0)[cc]{the tree $(T_1\hat{\ \ }T_2)$}}
\end{picture}
\vspace{1ex}

Picture~\theppp.
\end{center}

It follows from standard combinatorics that for any $n\ge0$, the number of rooted binary trees with $n$ carets is equal to the $n$th Catalan number $c_n=\frac{(2n)!}{n!(n+1)!}$.
\vspace{1ex}

Each rooted binary tree has {\em leaves\/}. Formally, they are defined as follows: for the trivial tree, the only leaf coincides with the root. In case $T=(T_1\hat{\ \ }T_2)$, the set of leaves equals the union of the sets of leaves of $T_1$ and $T_2$. In this case the leaves are exactly vertices of degree
$1$.

The following concept will play an important r\^ole throughout of the paper.

\begin{df}
\label{height}
The {\em height\/} of a rooted binary tree $T$ denoted by $\hgt T$ is defined by induction on the number of carets in $T$. For the trivial tree, its height equals $0$. For $T=(T_1\hat{\ \ }T_2)$, its height is $\hgt T=\max(\hgt T_1,\hgt T_2)+1$.
\end{df}

Now we define a {\em rooted binary forest\/} as a finite sequence of rooted binary trees $T_1$, ... , $T_m$,
where $m\ge1$. The leaves of it are the leaves of the trees. It is standard from combinatorics that the number
of rooted binary forests with $n$ leaves also equals $c_n$. The trees are enumerated from left to right and they
are drawn in the same way.

A {\em marked\/} (rooted binary) forest is a (rooted binary) forest where one of the trees is marked. On a picture we will usually show the marker as a vertical arrow on a tree.

\section{Semigroup diagrams and binary forests}
\label{sgpd}

Let ${\cal P}$ be a semigroup presentation. For our needs it suffices to work with the simplest presentation ${\cal P}=\la x\mid x=x^2\ra$. We will use 
non-spherical diagrams over ${\cal P}$ to represent elements of the group $F$. A more detailed description can be found in~\cite{Gu04}. Here we will mostly use diagrams that represent elements of the monoid $M$. Such objects are called {\em positive diagrams} over $x=x^2$. Let us give a brief illustration.

Given a normal form $pq^{-1}$ of an element in $F$, it is easy to draw the corresponding non-spherical diagram,
and vice versa. The following example illustrates the diagram that corresponds to the element $g=x_0^2x_1x_6x_3^{-1}x_0^{-2}$ represented by its normal form:

\refstepcounter{ppp}
\begin{center}
	\begin{picture}(87.00,30.00)
	\put(6.00,9.00){\circle*{1.00}}
	\put(16.00,9.00){\circle*{1.00}}
	\put(26.00,9.00){\circle*{1.00}}
	\put(36.00,9.00){\circle*{1.00}}
	\put(46.00,9.00){\circle*{1.00}}
	\put(56.00,9.00){\circle*{1.00}}
	\put(66.00,9.00){\circle*{1.00}}
	\put(76.00,9.00){\circle*{1.00}}
	\put(86.00,9.00){\circle*{1.00}}
	\put(6.00,9.00){\line(1,0){80.00}}
	\bezier{132}(66.00,9.00)(76.00,22.00)(86.00,9.00)
	\bezier{120}(16.00,9.00)(26.00,20.00)(36.00,9.00)
	\bezier{176}(36.00,9.00)(26.00,25.00)(6.00,9.00)
	\bezier{264}(46.00,9.00)(28.00,35.00)(6.00,9.00)
	\bezier{120}(36.00,9.00)(46.00,-2.00)(56.00,9.00)
	\bezier{104}(6.00,9.00)(16.00,1.00)(26.00,9.00)
	\bezier{176}(6.00,9.00)(20.00,-7.00)(36.00,9.00)
	\put(36.00,15.00){\makebox(0,0)[cc]{$x_0$}}
	\put(14.00,12.00){\makebox(0,0)[cc]{$x_0$}}
	\put(26.00,12.00){\makebox(0,0)[cc]{$x_1$}}
	\put(76.00,12.00){\makebox(0,0)[cc]{$x_6$}}
	\put(25.00,5.00){\makebox(0,0)[cc]{$x_0^{-1}$}}
	\put(48.00,6.50){\makebox(0,0)[cc]{$x_3^{-1}$}}
	\put(18.00,7.00){\makebox(0,0)[cc]{$x_0^{-1}$}}
	\end{picture}
\vspace{1ex}

Picture~\theppp.
\end{center}

Here we assume that each edge of the diagram is labelled by a letter $x$. These labels are usually not shown. The horizontal path in the picture cuts the diagram into two parts, {\em positive} and {\em negative}. The positive part represents the element $x_0^2x_1x_6\in M$. The top path of this positive diagram has label $x^4$, the bottom path of it has label $x^8$. So we have a positive $(x^4,x^8)$-diagram over ${\cal P}$. 

Given a positive $(x^m,x^n)$-diagram $\Delta_1$ and a positive $(x^k,x^l)$-diagram $\Delta_2$, one can {\em concatenate} them in the following natural way. If $n\le k$, then we identify the bottom path of $\Delta_1$ with the initial segment of the top path of $\Delta_2$ of the same length. If $n>k$, then the initial segment of the bottom path of $\Delta_1$ of length $k$ is identified with the top path of $\Delta_2$. The result of this operation is denoted by $\Delta_1\circ\Delta_2$ or simply as $\Delta_1\Delta_2$. It is well known that this operation is associative.

The operation of {\em concatenation} described above requires cancelling dipoles (mirror images of cells) in case of representing elements of the group $F$. However, working with elements of $M$, we do not need this additional operation at all.
\vspace{2ex}

Now let us show how positive diagrams can be replaced by marked forests when representing elements of $M$. We restrict ourselves to the set of positive $(x^m,x^{n+1})$-diagrams where $m\ge2$ and $n$ is a fixed number. Let $\Delta$ be such a diagram. We mark the second edge from the left of the top path of $\Delta$. Then we remove all edges starting at the leftmost vertex of $\Delta$. The result will be a positive diagram $\Delta'$, where one of the top edges is marked, and the bottom path has length $n$. Let us illustrate this by the following expample.

Let $w=x_0^3x_4x_5x_8x_{10}^2$ be an element of $M$. It is represented by the following positive diagram $\Delta$:

\refstepcounter{ppp}
\begin{center}
\begin{picture}(106.75,20.75)(0,0)
\multiput(2,2)(13.09375,0){8}{\line(1,0){13.09375}}
\put(2,2){\circle*{1.}}
\put(10,2){\circle*{1.}}
\put(18,2){\circle*{1.}}
\put(26,2){\circle*{1.}}
\put(34,2){\circle*{1.}}
\put(42,2){\circle*{1.}}
\put(50,2){\circle*{1.}}
\put(58,2){\circle*{1.}}
\put(66,2){\circle*{1.}}
\put(74,2){\circle*{1.}}
\put(82,2){\circle*{1.}}
\put(90,2){\circle*{1.}}
\put(98,2){\circle*{1.}}
\put(106,2){\circle*{1.}}
\qbezier(2,2)(11.125,16.375)(18,2)
\qbezier(2,2)(10,25.5)(25.75,2)
\qbezier(2,2)(18,44.75)(34,2)
\qbezier(42,2)(51.125,16)(58,2)
\qbezier(58,2)(49.375,29.875)(34,2)
\qbezier(66,2)(74.75,19)(82,2)
\qbezier(81.5,2)(91.75,11.375)(98,2)
\qbezier(81.75,2)(92.625,23.5)(106,2)
\end{picture}
\vspace{1ex}

Picture~\theppp.
\end{center}

After the operations of deleting edges and marking an edge on the top, we get the following result:

\refstepcounter{ppp}
\begin{center}
\begin{picture}(106.75,20.75)(0,0)
\multiput(10,2)(12.09375,0){8}{\line(1,0){12.09375}}
\put(10,2){\circle*{1.}}
\put(18,2){\circle*{1.}}
\put(26,2){\circle*{1.}}
\put(34,2){\circle*{1.}}
\put(42,2){\circle*{1.}}
\put(50,2){\circle*{1.}}
\put(58,2){\circle*{1.}}
\put(66,2){\circle*{1.}}
\put(74,2){\circle*{1.}}
\put(82,2){\circle*{1.}}
\put(90,2){\circle*{1.}}
\put(98,2){\circle*{1.}}
\put(106,2){\circle*{1.}}
\qbezier(42,2)(51.125,16)(58,2)
\qbezier(58,2)(49.375,29.875)(34,2)
\qbezier(66,2)(74.75,19)(82,2)
\qbezier(82,2)(91.75,11.375)(98,2)
\qbezier(82,2)(92.625,23.5)(106,2)
\put(48.5,24.75){\vector(0,-1){5.5}}
\end{picture}
\vspace{1ex}

Picture~\theppp.
\end{center}

Now we introduce a standard operation to switch between positive diagrams and rooted binary forests. Given a positive diagram, we take the set of midpoints of all its edges. For any positive $(x,x^2)$-cell we connect the midpoint of its top edge with the midpoints of two its bottom edges. The two connecting intervals form a caret. After we erase the original diagram, we get to a sequence of rooted binary trees. Notice that to any edge that was a part of both the top path and the bottom path of the positive diagram, we get a trivial tree consisting of a single vertex.

Here is the result of the above operation applied to $\Delta'$ from Pic.~\theppp:

\refstepcounter{ppp}
\begin{center}
\begin{picture}(102.5,23.75)(0,0)
\put(14,2){\circle*{1.}}
\put(22,2){\circle*{1.}}
\put(30,2){\circle*{1.}}
\put(38,2){\circle*{1.}}
\put(46,2){\circle*{1.}}
\put(54,2){\circle*{1.}}
\put(62,2){\circle*{1.}}
\put(70,2){\circle*{1.}}
\put(78,2){\circle*{1.}}
\put(86,2){\circle*{1.}}
\put(94,2){\circle*{1.}}
\put(102.25,2){\circle*{1.}}
\put(46,2.25){\line(0,1){.5}}
\multiput(46,2)(.033687943,.046099291){141}{\line(0,1){.046099291}}
\multiput(50.75,8.5)(.03370787,-.07022472){89}{\line(0,-1){.07022472}}
\multiput(37.75,2.5)(.0336879433,.0487588652){282}{\line(0,1){.0487588652}}
\multiput(47.25,16.25)(.033653846,-.076923077){104}{\line(0,-1){.076923077}}
\multiput(69.75,2.5)(.03358209,.041044776){134}{\line(0,1){.041044776}}
\multiput(74.25,8)(.033482143,-.051339286){112}{\line(0,-1){.051339286}}
\multiput(86,2.25)(.033730159,.045634921){126}{\line(0,1){.045634921}}
\multiput(90.25,8)(.033653846,-.052884615){104}{\line(0,-1){.052884615}}
\multiput(90.5,8)(.033536585,.045731707){164}{\line(0,1){.045731707}}
\multiput(96,15.5)(.033707865,-.074438202){178}{\line(0,-1){.074438202}}
\put(47.25,23.75){\vector(0,-1){5.25}}
\end{picture}
\vspace{1ex}

Picture~\theppp.
\end{center}	

One of the trees here is marked so we get a marked rooted binary forest. It is easy to apply inverse operations. Say, the marked forest from the following picture

\refstepcounter{ppp}
\begin{center}
\begin{picture}(134.25,21.5)(0,0)
\put(14,2){\circle*{1.}}
\put(22,2){\circle*{1.}}
\put(30,2){\circle*{1.}}
\put(38,2){\circle*{1.}}
\put(46,2){\circle*{1.}}
\put(54,2){\circle*{1.}}
\put(62,2){\circle*{1.}}
\put(70,2){\circle*{1.}}
\put(78,2){\circle*{1.}}
\put(86,2){\circle*{1.}}
\put(94,2){\circle*{1.}}
\put(102.25,2){\circle*{1.}}
\put(110,2){\circle*{1.}}
\put(118,2){\circle*{1.}}
\put(126,2.25){\circle*{1.}}
\put(134,2){\circle*{1.}}
\multiput(14,2.25)(.03358209,.042910448){134}{\line(0,1){.042910448}}
\multiput(18.5,8)(.033653846,-.057692308){104}{\line(0,-1){.057692308}}
\multiput(38,2.25)(.033557047,.038590604){149}{\line(0,1){.038590604}}
\multiput(43,8)(.03370787,-.06460674){89}{\line(0,-1){.06460674}}
\put(54.25,2.5){\line(2,3){3.5}}
\multiput(57.75,7.75)(.033613445,-.046218487){119}{\line(0,-1){.046218487}}
\multiput(43,7.75)(.033613445,.038865546){238}{\line(0,1){.038865546}}
\multiput(51,17)(.033678756,-.046632124){193}{\line(0,-1){.046632124}}
\multiput(86.25,2.5)(.033730159,.043650794){126}{\line(0,1){.043650794}}
\multiput(90.5,8)(.033505155,-.059278351){97}{\line(0,-1){.059278351}}
\multiput(93.75,2.25)(-.03370787,.06179775){89}{\line(0,1){.06179775}}
\put(90.75,7.75){\line(3,4){4.5}}
\multiput(95.25,13.75)(.033678756,-.056994819){193}{\line(0,-1){.056994819}}
\multiput(125.75,2.25)(.033653846,.040064103){156}{\line(0,1){.040064103}}
\multiput(131,8.5)(.033505155,-.069587629){97}{\line(0,-1){.069587629}}
\multiput(117.75,2)(.0336538462,.0392628205){312}{\line(0,1){.0392628205}}
\multiput(128.25,14.25)(.03370787,-.07022472){89}{\line(0,-1){.07022472}}
\multiput(110,2)(.0336842105,.0389473684){475}{\line(0,1){.0389473684}}
\multiput(126,20.5)(.03358209,-.09328358){67}{\line(0,-1){.09328358}}
\put(95.25,21.5){\vector(0,-1){6}}
\end{picture}

\vspace{1ex}

Picture~\theppp.
\end{center}

\noindent will lead to the element $w=x_0^5x_1x_4^2x_6x_{10}^2x_{13}x_{14}x_{15}\in M$.
\vspace{1ex}

Any set of marked binary forests with $n$ leaves is a finite subset in the Cayley graph of $F$. We prefer to work with left Cayley graphs since in this case we need to consider positive elements instead of negative ones. Otherwise we would work with mirror images of rooted binary forests with respect to a horizontal axis. Recall that an edge labelled by $a\in A^{\pm1}$, where $A$ is the set of group generators, starting at a vertex $g$, goes to the element $a^{-1}g$ in the left Cayley graph. We need an explicit description of acting all of the generators $x_0^{\pm1}$, $x_1^{\pm1}$ for the case of the left Cayley graph of $F$ in these generators.

Here are the rules of the game.

Let $g$ be an element represented by a marked rooted binary forest $\dots,T_{-1},T_0,T_1,\dots$, where $T_0$ is the marked tree. 

\begin{itemize}
\item
Acting by $x_0$, that is, going along the edge with label $x_0$, means to move the marker left. That is, $T_{-1}$ becomes the marked tree for the same forest. This operation is not applied if $T_0$ is the leftmost tree.

\item
Acting by $x_0^{-1}$, that is, going along the edge with label $x_0^{-1}$, means to move the marker right. That is, $T_1$ becomes the marked tree for the same forest. This operation is not applied if $T_0$ is the rightmost tree.

\item
Acting by $x_1$, that is, going along the edge with label $x_1$, means to remove the caret from the tree $T_0=(T'\hat{\ \ }T'')$. The marker will have its position at $T'$. This operation is not applied if $T_0$ is a trivial tree.

\item
Acting by $x_1^{-1}$, that is, going along the edge with label $x_1^{-1}$, means to put a new caret over $T_0$ and $T_1$ making this tree marked. That is, instead of $T_0$ and $T_1$ in the forest we get the tree $T=(T_0\hat{\ \ }T_1)$. This operation is not applied if $T_0$ is the rightmost tree.
\end{itemize}

This description will be important for constructing finite subgraphs in the left Cayley graph of $F$.

\section{Amenability and the density of the Cayley graphs}
\label{dens}

Recall that the group is called {\em amenable} whenever there exists a finitely additive normalized invariant mean on $G$, that is, a mapping $\mu\colon{\cal P}(G)\to[0,1]$ such that
\begin{itemize}
\item
 $\mu(Z_1\cup Z_2)=\mu(Z_1)+\mu(Z_2)$ for any disjoint subsets $Z_1,Z_2\subseteq G$ 
\item 
$\mu(G)=1$
\item
$\mu(Zg)=\mu(gZ)=\mu(Z)$ for any $Z\subseteq G$, $g\in G$.
\end{itemize}

One gets an equivalent definition of amenability if only one-sided invariance of the mean is assumed, say, the condition $\mu(Zg)=\mu(Z)$ ($Z\subseteq G$, $g\in G)$. The proof can be found in \cite{GrL}.

We are not going to list all well-known properties of (non)amenable groups. It is sufficient to refer to one of modern surveys like~\cite{Sap14}. Just notice that all finite and all Abelian groups are amenable. The class of amenable groups is closed under taking subgroups, homomorphic images, group extensions, and directed unions of groups. The groups in the closure of the classes of the union of finite and Abelian groups under this list of operations are called {\em elementary amenable} (EA). Also we recall that free groups of rank $> 1$ are not amenable. 
\vspace{1ex}

We will often refer to the following well-known statement equivalent to the F\o{}lner criterion~\cite{Fol}. Here we restrict ourselves to the case of finitely generated groups.

\begin{prop}
\label{fol}
A group $G$ with finite set of generators $A$ is amenable if and only if its Cheeger isoperimetric constant is zero: $\iota_*(G;A)=0$. Equivalently, if $|A|=m$, then the Cayley graph ${\cal C}(G;A)$ has density $2m$.
\end{prop}

This holds for any finite set of generators and can be applied to both Cayley graphs: right or left.

In practice, to establish amenability of a group, it is sufficient to construct a collection of finite subgraphs $Y$ in the Cayley graph such that $\inf\limits_Y\frac{|\partial_{\ast}Y|}{|Y|}=0$. Such subsets of vertices are called {\em Folner sets}. Informally, this means that almost all vertices of these sets are internal.
\vspace{1ex}

Brin and Squier~\cite{BS} proved that the group $F$ has no free subgroups of rank $>1$. It is also known that $F$ is not elementary amenable~\cite{Chou}. However, the famous problem about amenability of $F$ is still open. The question whether $F$ is amenable was asked by Ross Geoghegan in 1979; see~\cite{Geo,Ger87}. There were many attempts of various authors to solve this problem in both directions. Some of these papers contained new interesting ideas. We are not going to review a detailed history of the problem; this information can be found elsewhere. However, to emphasize the difficulty of the question, we mention the paper~\cite{Moore13}, where it was shown that if $F$ is amenable, then Folner sets for it have a very fast growth (like towers of exponents). 

Notice that if $F$ is amenable, then it brings an example of a finitely presented amenable group which is not EA. If it is non-amenable, then this gives an example of a finitely presented group, which is not amenable and has no free subgroups of rank $>1$.  Notice that the first example of a non-amenable group without free non-abelian subgroups has been constructed by Ol'shanskii \cite{Olsh}. (The question about such groups was formulated in \cite{Day}, it is also often attributed to von Neumann \cite{vNeu}.) Adian \cite{Ad83} proved that free Burnside groups with $m>1$ generators of odd exponent $n\ge665$ are not amenable. The first example of a finitely presented non-amenable group without free non-abelian subgroups has been constructed by Ol'shanskii and Sapir \cite{OlSa}. Grigorchuk \cite{Gri} constructed the first example of a finitely presented amenable group which is not EA.
\vspace{1ex}

Now let us discuss the question about the density of the Cayley graph of $F$ in the standard generating set $A=\{x_0,x_1\}$. It was shown in~\cite{Gu04} that the density is at least $3$. In the Addendum to the same paper, there was a modification of the above construction showing that there are finite subgraphs with density strictly greater than $3$. An essential improvement was obtained in~\cite{BB05}. The authors constructed a family of finite subgraphs with density approaching $3{.}5$. This was the best known estimate for the density of the Cayley graph of $F$ for many years. Some authors believed that the above estimate was exact (for instance, see~\cite{Bur16}). There were several heuristic reasons for that. See~\cite{Gu22} for the corresponding Conjecture 1. However, it turned out that $3{.}5$ is not the best estimate for the density. One of the main results of the present paper can be stated as follows.

\begin{thm}
\label{over3.5}
The density of the Cayley graph of R.\,Thompson's group $F$ in the standard set of generators $\{x_0,x_1\}$ is strictly greater than $3{.}5$. Equivalenly, the Cheeger isoperimetric constant of the group $F$ in the same set of generators is stricltly less than $\frac12$.
\end{thm}

This result increases the chances that $F$ is amenable.
\vspace{1ex}

We will need a precise description of the Belk -- Brown construction. The reader can find more information in~\cite{Be04,BB05}. Notice that in~\cite{Gu21a} we also gave a description of Belk -- Brown sets.

The following elementary lemma was proved in~\cite{Gu21a}. We will refer to it as the {\em symmetric property}.

\begin{prop}
\label{invlet}
Let $G$ be a finitely generated group and let $\Gamma={\cal C}(G,A)$ be its Cayley graph. Let $Y$ be a finite nonempty subgraph of $\Gamma$. Then for any $a\in A^{\pm1}$, the number of edges in the Cheeger boundary $\partial_{\ast} Y$ labelled by $a$ is the same as the number of edges in $\partial_{\ast}Y$ labelled by $a^{-1}$.
\end{prop}

Let $n\ge1$, $k\ge0$ be integer parameters. By $BB(n,k)$ we denote the set of all marked forests that have $n$ leaves, and
each of their trees has height at most $k$. According to Section~\ref{sgpd}, we regard $BB(n,k)$ as a finite subset of $F$, that is, a set of vertices of the left Cayley graph of $F$ in standard generators. These generators have a partial action on this set. The rules are described at the end of the previous Section. We only remark that when we add a caret over two trees (acting by $x_1^{-1}$), we need to claim that each of these trees have height at most $k-1$. Otherwise we go outside the set $BB(n,k)$ going along the edge with label $x_1^{-1}$. 

Let us discuss briefy some properties of these sets including heuristic reasons for the conjecture that the above construction could be optimal. For any fixed $k$, let $n\gg k$. Since any tree of height $k$ has at most $2^k$ leaves, any forest in $BB(n,k)$
contains at least $\frac{n}{2^k}$ trees. So if we choose a marked forest in random, the probabililty for
this vertex of an automaton to accept both $x_0$, $x_0^{-1}$ approaches $1$. Now look at the probability to
accept $x_1$. The contrary holds if and only if the marked tree is trivial. When we remove this trivial tree, we just obtain an element of $BB(n-1,k)$. So the probability we are interested in, is the quotient $\frac{|BB(n-1,k)|}{|BB(n,k)|}$. But $BB(n,k)$ is the set of forests with $n$ leaves. The number of them is the Catalan number $c_n\le4^n$. So one cannot expect the value of the quotient better than $\frac14$. 
The probability not to accept the inverse letter $x_1^{-1}$ is the same by the symmetric property (Proposition~\ref{invlet}). This gives an expected bound $\frac12$ for the Cheeger isoperimetric constant. So many people believed that the value $3.5$ for the density of the Cayley graph of $F$ should be optimal. 
\vspace{1ex}

One possible idea to modify the construction was to change the lower bound for the height of trees in the forest. According to our private communication with Jim Belk, this approach does not work. It gives either the same bound or even worse.

However, the above sets can be modified in a different way using some ideas of~\cite{AGL08} and also our recent papers~\cite{Gu21a,Gu22}. The details will be given in the next Sections.

\section{Generating Functions}
\label{gf}

Here we are giving details on the generating functions of the sets of forests and marked forests. We need this to obtain analytic estimates for the size of some finite subsets in the Cayley graph of $F$. General information about the techniques of generating functions can be found in~\cite{Wilf}.
\vspace{1ex}

Let us define a sequence of polynomials by induction:
\be{phi0}
\Phi_0(x)=x
\ee
\be{phik}
\Phi_k(x)=x+\Phi_{k-1}(x)^2\quad\quad k\ge1.
\ee

Notice that $\Phi_k(x)$ is the generating polynomial for the set of trees of height at most $k$. This means that the coefficient on $x^n$ in this polynomial shows the number of such trees with $n$ leaves. This follows directly from~(\ref{phik}). The summand $x$ corresponds to the trivial tree (with one leaf); for the tree $T=(T_1\hat{\ \ }T_2)$ of height $\le k$ we have height $\le k-1$ for each of the trees $T_1$, $T_2$. By induction, the pair of them has generating function $\Phi_{k-1}(x)^2$. This agrees with~(\ref{phik}).
\vspace{1ex}

It is easy to see that the equation $\Phi_k(x)=1$ has a unique positive root that we denote by $\xi_k$. Let us give an estimate for it. First of all, $\Phi_k(\frac14) < \frac12$ by induction. Thus $\frac14 < \xi_k$. On the other hand, $\Phi_k(\frac14)\ge\frac12-\frac1{k+4}$, where the inequality is strict whenever $k\ge1$. Indeed, this holds for $k=0$, and for $k\ge1$ one has by induction
$$
\Phi_k\left(\frac14\right)\ge\frac14+\left(\frac12-\frac1{k+3}\right)^2=\frac12-\frac{k+2}{(k+3)^2} > \frac12-\frac1{k+4}\,.
$$

Let $\varepsilon>0$. We show that $\Phi_k(\frac14+\varepsilon)\ge\Phi_k(\frac14)+w_k\varepsilon$, where $w_k=\frac{k+4}3-\frac2{(k+2)(k+3)}$.  This holds for $k=0$ since $w_0=1$. For $k\ge1$ we have by induction on $k$:
$$
\Phi_k\left(\frac14+\varepsilon\right)\ge\frac14+\varepsilon+\left(\Phi_{k-1}\left(\frac14\right)+w_{k-1}\varepsilon\right)^2 > \frac14+\Phi_{k-1}\left(\frac14\right)^2+\varepsilon\left(1+2\Phi_{k-1}\left(\frac14\right)w_{k-1}\right).
$$

Using that $2\Phi_{k-1}(\frac14)\ge1-\frac2{k+3}$, we see that the coefficient on $\varepsilon$ is greater than or equal to $1+\frac{k+1}{k+3}(\frac{k+3}3-\frac2{(k+1)(k+2)})=w_k$, so $\Phi_k(\frac14+\varepsilon)>\Phi_k(\frac14)+w_k\varepsilon$, as desired.
\vspace{1ex}

Taking $\varepsilon=\frac3{2k}$, one can check that $\Phi_k(\frac14+\varepsilon)>\frac12-\frac1{k+4}+\frac{3w_k}{2k} > 1$. This implies an upper bound: $\xi_k < \frac14+\frac3{2k}$. In particular, $\xi_k\to\frac14$ as $k\to\infty$.

Notice that for $k\gg1$ the value of $\xi_k$ is close to $\frac14+\frac{\pi^2}{k^2}$ as some numerical experiments show. However, we do not need any exact estimates here.
\vspace{1ex}

Let $R$ be the set of all complex roots of the equation $\Phi_k(x)=1$ different from $\xi_k$. By the triangle inequality, $1=|\Phi_k(\lambda)|\le\Phi_k(|\lambda|)$, where the inequality is strict for all $\lambda\in R$ since no other roots are positive reals. Hence $|\lambda|>\xi_k$ for all $\lambda\in R$. It is also clear that $\xi_k$ is a simple root of $\Phi_k(x)$ since the polynomial is increasing on $x>0$.

It is easy to see that the degree of $\Phi_k$ equals $m=2^k$. Let
$$
\Phi_k(x)=p_1x+p_2x^2+\cdots+p_mx^m.
$$
It is clear that all coefficients $p_i$ ($1\le i\le m$) are strictly positive, where $p_1=p_m=1$. We are interested in the following generating function
\be{1phik}
\frac1{1-\Phi_k(x)}=1+\Phi_k(x)+\Phi_k(x)^2+\cdots=\sum\limits_{n\ge0}\alpha_k(n)x^n
\ee
presented as a power series. Clearly, $\Phi_k(x)^l$ is the generating function for the set of forests with $l\ge0$ trees, where all these trees have height at most $k$. Adding all these polynomials, we get~(\ref{1phik}). So the coefficients $\alpha_k(n)$ show the number of forests (not marked) having $n$ leaves with the upper bound $k$ for the height of each of their trees. Obviously, all these numbers are strictly positive. These coefficients satisfy the following linear recurrent equations
$$
\alpha_k(n+m)=p_1\alpha_k(n+m-1)+\cdots+p_m\alpha_k(n)
$$
for $n\ge0$. It follows from the theory of such equations~\cite[Theorem 2.2]{SeFl} that the space of their solutions has a basis formed by sequences $\mu^n$, $n\mu^n$, \dots, $n^{d-1}\mu^n$, where $\mu$ runs over the set of roots of the characteristic equation for the above recurrence, where $d$ is the multiplicity of the root $\mu$. In our case the characteristic equation is
$$
\mu^{m}-p_1\mu^{m-1}-\cdots-p_{m-1}\mu-p_m=0
$$
so its roots are $\xi_k^{-1}$ and $\lambda^{-1}$ ($\lambda\in R$). The theory allows us to write down $$
\alpha_k(n)=C(\xi_k^{-1})^n+\sum\limits_{\lambda\in R}Q_{\lambda}(n)(\lambda^{-1})^n
$$
for some polynomials $Q_{\lambda}(n)$. We already know that $|\xi_k^{-1}|>|\lambda^{-1}|$ for all $\lambda\in R$. Let us show that $C > 0$. 

The coefficients $\alpha_k(0)$, $\alpha_k(1)$, \dots, $\alpha_k(m-1)$ are positive so there exits a positive constant $\rho$ such that $\alpha_k(n)\ge\rho(\xi_k^{-1})^n$ for $n\in\{0,1,\ldots,m-1\}$. Applying the recurrent equation and using that $\xi_k^{-1}$ satisfies the characteristic equation, we see that $\alpha_k(n)\ge\rho(\xi_k^{-1})^n$ for all $n\ge0$ by induction. Therefore, $C$ is nonzero since otherwise $\alpha_k(n)=o((\xi_k^{-1})^n)$. So $C$ equals the limit $\alpha_k(n)\xi_k^n$ as $n\to\infty$ and therefore it is a positive real number.

We obtained that $\alpha_k(n)\sim C(\xi_k^{-1})^n$ where the equavalence means that the quotient is asymptotically $1$. In particular, $\lim\limits_{n\to\infty}\frac{\alpha_k(n-1)}{\alpha_k(n)}=\xi_k$. The existence of such a limit is not true in general if we replace $\Phi_k(x)$ by arbitrary polynomials.
\vspace{1ex}

Another important generating function will be
\be{phik2}
\frac{\Phi_k(x)}{(1-\Phi_k(x))^2}=\sum\limits_{n\ge0}\beta_k(n)x^n.
\ee

Here $\Phi_k(x)$ is the generating fuction for the choice of a tree that will be marked. The generating function for forests to the right from the marked tree, as well as the forests to the left of it, equals~(\ref{1phik}). So~(\ref{phik2}) gives us the generating function for the set of marked forests with the bound $k$ for the height. The coefficient $\beta_k(n)$ shows the number of such forests with $n$ leaves. That is, this is a cardinality of the Belk -- Brown set: $|BB(n,k)|=\beta_k(n)$. 

To find the growth rate of the coefficients of this and other power series, we need the following

\begin{lm}
\label{grrs}
Let $P(x)$ be a polynomial, and let $\Psi(x)=\sum\limits_{n\ge0}\gamma(n)x^n$ be a power series, where all coefficients are positive, and the limit $\xi=\lim\limits_{n\to\infty}\frac{\gamma(n-1)}{\gamma(n)}$ exists. Then the power series for the product $P(x)\Psi(x)=\sum\limits_{n\ge0}\delta(n)x^n$ satisfies the following condition: $$\lim\limits_{n\to\infty}\frac{\delta(n)}{\gamma(n)}=P(\xi).$$
\end{lm}

\prf Let $P(x)=a_0+a_1x+\cdots+a_sx_s$. Comparing the coefficients on $x^n$ for the formal series, we get $\delta(n)=a_0\gamma(n)+a_1\gamma(n-1)+\cdots+a_s\gamma(n-s)$. Dividing by $\gamma(n)$ and using the fact that $\lim\limits_{n\to\infty}\frac{\gamma(n-r)}{\gamma(n)}=\xi^r$ for any $r\ge1$, we obtain that $\frac{\delta(n)}{\gamma(n)}\to a_0+a_1\xi+\cdots+a_s\xi^s=P(\xi)$ as $n\to\infty$.

The proof is complete.

Although the fact is elementary, we will widely use this lemma throughout the paper in the process of estimating coefficients of generating functions.
\vspace{1ex}

{\bf Remark.} Lemma~\ref{grrs} is very useful to estimate probabilities of some events. We can be interested to find the probability that the marked tree belongs to some finite set. For instance, that the marked tree has height $k$. The generating polynomial for the set of these trees is $P(x)=\Phi_k(x)-\Phi_{k-1}(x)$. If we take the generating function
\be{1phik2}
\frac1{(1-\Phi_k(x))^2}=\sum\limits_{n\ge0}\gamma_k(n)x^n
\ee
and multiply it by $P(x)$, then we get the generating function for the set of forests in which the tree of height $k$ is marked. According to Lemma~\ref{grrs}, the limit of the above probability as $n\to\infty$ will be $P(\xi_k)=1-\sqrt{1-\xi_k}$. The same concerns some other probabilities. Say, in~\cite{Gu21a} we studied the Cayley graph of $F$ in symmetric generators $\{x_1,\bar{x}_1\}$, where $\bar{x}_1=x_1x_0^{-1}$. We needed to estimate the probability that the marked tree has height $k$ together with the tree to the right of it. We did it in a combinatorial way. Now the same can be done analitically taking the square of the above polynomial $P(x)$ and the value of it at $x=\xi_k$. In general, Lemma~\ref{grrs} shows that separate events of the above form are asymtotically independent, that is, we can simply multiply the limits of probabilities for them.
\vspace{1ex}

We will be interested in finding the asymptotics of the coefficients of~(\ref{1phik2}). It is easy to see that the characteristic equation for the linear recurrent relations of $\gamma_k(n)$ is the square of the equation for $\alpha_k(n)$, that is,
$$
(\mu^{m}-p_1\mu^{m-1}-\cdots-p_{m-1}\mu-p_m)^2=0.
$$
Therefore, the complex roots of it will be the same, with double multiplicity. In particular, $\xi_k^{-1}$ will have multiplicity two. According to~\cite[Theorem 2.2]{SeFl}, we have 
$$\gamma_k(n)=(An+B)(\xi_k^{-1})^n+\sum\limits_{\lambda\in R}S_{\lambda}(n)(\lambda^{-1})^n$$
for some polynomials $S_{\lambda}(n)$. We remember that $|\lambda^{-1}| < \xi_k^{-1}$ for all $\lambda\in R$. On the other hand, we know that $\alpha_k(n)\ge\rho(\xi_k^{-1} )^n$ for some positive constant $\rho$ for all $n\ge0$. Since~(\ref{1phik2}) is the square of~(\ref{1phik}), one has
$$
\gamma_k(n)=\alpha_k(0)\alpha_k(n)+\alpha_k(1)\alpha_k(n-1)+\cdots+\alpha_k(n)\alpha_k(0)\ge(n+1)\rho^2(\xi_k^{-1})^n.
$$
This implies $A\ne0$. Taking the limit of $\frac{\gamma_k(n)\xi_k^n}n$ as $n\to\infty$, we see that $A$ is a positive real number. So $\gamma_k(n)\sim An(\xi^{-1})^n$ and the limit $\lim\limits_{n\to\infty}\frac{\gamma_k(n-1)}{\gamma_k(n)}=\xi_k$ exists. Now we can apply Lemma~\ref{grrs} taking~(\ref{1phik2}) for $\Psi(x)$ and $\Phi_k(x)$ for $P(x)$. This gives us 
$\lim\limits_{n\to\infty}\frac{\beta_k(n)}{\gamma_k(n)}=\Phi_k(\xi_k)=1$. Therefore, $\beta_k(n)\sim\gamma_k(n)\sim An(\xi_k^{-1})^n$, and we know the asymptotic growth of $\beta_k(n)$. In particular, $\lim\limits_{n\to\infty}\frac{\beta_k(n-1)}{\beta_k(n)}=\xi_k$.

\section{Proof of Theorem 1}
\label{thm1}

The proof will essentially use almost all notation introduced in the previous Section. The basic idea of the proof is as follows. The subset $BB(n,k)$ of the left Cayley graph of $F$ has some fragments with small density. The probability to meet such a fragment has a positive uniform lower bound. If we remove such fragments from the subgraph, we increase its density exceeding the value $3{.} 5$. 

Let us describe the fragments we are interested in. Let we have a rooted binary forest $\dots,T_0,T_1,T_2,T_3,T_4,\dots$, where $T_0$ is marked and all trees have height at most $k$. We claim that trees $T_i$ ($0\le i\le4$) exist in this forest, and the following conditions hold:

\begin{itemize}
\item $T_0$ and $T_2$ are trivial trees,
\item $T_1$ and $T_3$ have height $k$,
\item $T_4$ is a nontrivial tree.
\end{itemize}

The latter condition is added for simplicity. A forest satisfying the listed conditions will be called {\em special}. To each of these forests we assign vertices $a$, $b$, $c$ of the left Cayley graph of $F$ in the standard generating set. Here $a$ corresponds to the forest with $T_0$ as a marked tree; $b$ and $c$ denote the forests where $T_1$ and $T_2$ are marked trees, respectively.

In the left Cayley graph of $F$ in standard generators these vertices look as follows:

\refstepcounter{ppp}
\begin{center}
\begin{picture}(87.75,18.552)(0,0)
\put(23.75,13.75){\circle{9.014}}
\put(46.25,13.75){\circle{9.014}}
\put(69,13.5){\circle{9.014}}
\put(42.25,13.75){\vector(-1,0){14.25}}
\put(64.75,13.75){\vector(-1,0){14.25}}
\put(19.5,13.75){\vector(-1,0){11.75}}
\put(87.75,13.75){\vector(-1,0){14.5}}
\put(46,9.5){\vector(0,-1){9.5}}
\put(23.5,13.5){\makebox(0,0)[cc]{$a$}}
\put(46.25,13.75){\makebox(0,0)[cc]{$b$}}
\put(69.,13.25){\makebox(0,0)[cc]{$c$}}
\put(13.75,16.75){\makebox(0,0)[cc]{$x_0$}}
\put(36,16){\makebox(0,0)[cc]{$x_0$}}
\put(57.75,15.5){\makebox(0,0)[cc]{$x_0$}}
\put(80.75,15.5){\makebox(0,0)[cc]{$x_0$}}
\put(50,3.25){\makebox(0,0)[cc]{$x_1$}}
\end{picture}
\vspace{1ex}

Picture~\theppp.
\end{center}

Since the tree $T_0$ is empty, we cannot remove a caret from it. This means that $a$ does not accept $x_1$ in the automaton $BB(n,k)$. The tree $T_1$ has height $k$ so a caret cannot be added to $T_0$ and $T_1$ to stay within $BB(n,k)$. So $a$ can accept only letters $x_0$ and $x_0^{-1}$. Exactly the same situation holds for the vertex $c$. Notice that the leftmost edge labelled by $x_0$ may not belong to the subgraph if $T_0$ is the leftmost tree in the special forest.

As for the vertex $b$, we can remove the caret from $T_1$ so $b$ accepts $x_1$. However, it does not accept $x_1^{-1}$ since the tree $T_1$ has height $k$ and no caret can be added to $T_1$ and $T_2$. Thus $a$, $c$ have degree 2 in $BB(n,k)$ and $b$ has degree 3 in the same subgraph.

If we have another special forest with the corresponding vertices $a'$, $b'$, $c'$, then no coincidences of vertices can occur. The only case could be $a'=c$ (or $a=c'$, which is totally symmetric). However, this is impossible by the choice of the tree $T_4$ in the special forest. If we go from $a'$ by a path labelled by $x_0^{-2}$, then we meet vertex $c'$ that corresponds to the trivial tree. Going along the path with the same label from $c$, we get the forest with marked tree $T_4$, which is nontrivial by definition. So this condition allows us to avoid repetitions.

Now let us find the number $\sigma_k(n)$ of special forests in $BB(n,k)$ using generating functions. The set consisting of one trivial tree has generating function $x$. This concerns trees $T_0$ and $T_2$. The generating functions for the set of nontrivial trees of height $\le k$ is $\Phi_k(x)-x$. This concerns $T_4$. The generating functions for the set of trees of height $k$ is $\Phi_k(x)-\Phi_{k-1}(x)$. This concerns trees $T_1$ and $T_3$. Taking a product 
$$P(x)=x^2(\Phi_k(x)-\Phi_{k-1}(x))^2(\Phi_k(x)-x),$$ 
we get a polynomial for the set of forests $T_0$, \dots, $T_4$ with the above conditions. Adding forests to the left and to the right of these 5 trees, we get the generating function $\frac{P(x)}{(1-\Phi_k(x))^2}$.

The number of special forests divided by the coefficient $\gamma_k(n)$ of~(\ref{1phik2}), approaches $P(\xi_k)$ according to Lemma~\ref{grrs}. That is, $\lim\limits_{n\to\infty}\frac{\sigma_k(n)}{\gamma_k(n)}=P(\xi_k)$. We know that $\Phi_{k-1}(x)=\sqrt{\Phi_k(x)-x}$ by~(\ref{phik}), so $P(\xi_k)=\xi_k^2(1-\sqrt{1-\xi_k})^2(1-\xi_k)$. When $k\to\infty$, this value approaches $p=\frac3{64}(1-\frac{\sqrt3}2)^2>\frac1{1200}$. This shows that special forests occur in $BB(n,k)$ with guaranteed positive probability that does not depend on the parameters whenever $n\gg k\gg1$.

Now let us exclude from $BB(n,k)$ triples of vertices $a$, $b$, $c$ for each special forest. Also we delete all directed edges incident to some of the deleted vertices. The result will be a subgraph in the left Cayley graph of $F$. We denote it by $BB'(n,k)$. Let us estimate its density. For every special forest, we delete 3 vertices and no more than 5 geometric edges, that is, no more than 10 directed edges. So the set of vertices of $BB'(n,k)$ has cardinality $|BB'(n,k)|=|BB(n,k)|-3\sigma_k(n)=\beta_k(n)-3\sigma_k(n)$.

The set of directed egdes of the whole Cayley graph incident to vertices in $BB(n,k)$ has cardinality $4\beta_k(n)$. Some of these edges are external, the others are internal. The number of external edges labelled by $x_0$ is $\alpha_k(n)$; they correspond to forests with the leftmost tree being marked. The generating function for them is~(\ref{1phik}). The same hold for external edges with label $x_0^{-1}$: they correspond to the forests with the rightmost tree being marked. So the number of external edges in $B(n,k)$ labelled by $x_0^{\pm1}$ equals $2\alpha_k(n)$. Recall that $\alpha_k(n)\sim C(\xi_k^{-1})^n$ and $\beta_k(n)\sim An(\xi_k^{-1})^n$ so $\alpha_k(n)=o(\beta_k(n))$. 

Vertices in $BB(n,k)$ that do not accept $x_1$ correspond to the case when the marked tree is trivial. The generating function for the set of these vertices is $\frac{x}{(1-\Phi_k(x))^2}$. So the coefficient on $x^n$ is $\beta_k(n-1)$. Exactly the same number of vertices do not accept $x_1^{-1}$ according to symmetric property (Proposition~\ref{invlet}). The same fact can be also checked directly using generating functions.

So we have $2\alpha_k(n)+2\beta_{k-1}(n)$ external edges in $BB(n,k)$. The other $4\beta_k(n)-2\alpha_k(n)-2\beta_{k-1}(n)$ edges are internal. Dividing by the number of vertices, we get $4-2\frac{\beta_k(n-1)}{\beta_k(n)}-2\frac{\alpha_k(n)}{\beta_k(n)}$ approaching $4-2\xi_k$ as $n\to\infty$ and then approaching $3{.}5$ as $k\to\infty$. This explains why the density of Belk -- Brown sets approach $3{.}5$.

As for our sets $BB'(n,k)$, they have $\beta_k(n)-3\sigma_k(n)$ vertices. The number of directed edges in $BB'(n,k)$ is greater than or equal to $4\beta_k(n)-2\alpha_k(n)-2\beta_{k-1}(n)-10\sigma_k(n)$. So the density of $BB'(n,k)$ has a lower bound
$$
\frac{4\beta_k(n)-2\alpha_k(n)-2\beta_{k-1}(n)-10\sigma_k(n)}{\beta_k(n)-3\sigma_k(n)}=
$$
$$\frac{4-2\frac{\alpha_k(n)}{\beta_k(n)}-2\frac{\beta_k(n-1)}{\beta_k(n)}-10\frac{\sigma_k(n)}{\beta_k(n)}}{1-3\frac{\sigma_k(n)}{\beta_k(n)}}\to\frac{4-2\xi_k-10P(\xi_k)}{1-3P(\xi_k)}
$$
as $n\to\infty$. The latter approaches 
$$
\frac{3{.}5-10p}{1-3p}=3{.5}+\frac{0{.}5p}{1-3p} > 3{.}5+\frac1{2400}
$$
as $k\to\infty$. So for some $n\gg k\gg1$ we can reach the density of finite subgraphs of the Cayley graph of $F$ in standard generators greater than $3{.}5004$.

This completes the proof of Theorem~\ref{over3.5}.
\vspace{1ex}

Notice that sets $BB'(n,k)$ for which the densities are close to this value, are really huge. Indeed, we had an estimate $\xi_k < \frac14+\frac3{2k}$ in Section~\ref{gf}. Our inequalities will work provided $0{.}5p > \frac3k$ so $k$ will be around $7200$. For the value of $n$ we need at least $n>2^k$ to get a large number of trees in a marked forest (or even higher). Finally, the number of forests with $n$ leaves grow like a Catalan number, that is, around $4^n$. Therefore, our proof works for sets of size greater than $2^{2^{7200}}$. Even if we improve the estimate of $\xi_k$, the exponent $7200$ will be replaced by several hundreds. It is clear that such sets are out of any reasonable computer search. These sizes look like towers of exponents from~\cite{Moore13}.

\section{Proof of Theorem 2}
\label{thm2}

The following result disproves our Conjecture 2 from~\cite{Gu22}.

\begin{thm}
\label{dens012}
The Cheeger isoperimetric constant of the Cayley graph of Thompson's group $F$ in the generating set $\{x_0,x_1,x_2\}$ is strictly less than $1$. Equivalently, the density of the corresponding Cayley graph strictly exceeds $5$.
\end{thm}

According to~\cite[Theorem 1]{Gu22}, this means that there are no pure evacuation schemes on the Cayley graph of $F$ in these generators.
\vspace{1ex}

The proof is based on the same idea as the one for Theorem~\ref{over3.5}. We take the sets of vertices $BB(n,k)$ and consider special forests defined in the beginning of Section~\ref{thm1}. Now the set of our generators is $\{x_0,x_1,x_2\}$ so we have to describe the conditions to accept $x_2$ and $x_2^{-1}$ for a vertex in an automaton. For a marked forest, accepting $x_2$ means that the tree to the right of the marker exists and it is nontrivial. Applying $x_2$ means that we remove the caret from this tree, and our marker stays at its original place.

To accept $x_2^{-1}$ for a vertex means that there exist two trees to the right of the marker, and each of them has height $\le k-1$. In this case we can add a caret over this pair of trees staying inside $BB(n,k)$. The marker does not change its place.

According to this description, the fragment of the left Cayley graph of $F$ in generators $\{x_0,x_1,x_2\}$ will look as follows for any special forest (we keep the notions and the notation of the previous Section):

\refstepcounter{ppp}
\begin{center}
\begin{picture}(87.75,18.553)(0,0)
\put(23.75,13.75){\circle{9.014}}
\put(46.25,13.75){\circle{9.014}}
\put(69,13.5){\circle{9.014}}
\put(42.25,13.75){\vector(-1,0){14.25}}
\put(64.75,13.75){\vector(-1,0){14.25}}
\put(19.5,14){\vector(-1,0){11.75}}
\put(87.75,13.75){\vector(-1,0){14.5}}
\put(46,9.5){\vector(0,-1){9.5}}
\put(23.5,13.5){\makebox(0,0)[cc]{$a$}}
\put(46.25,13.75){\makebox(0,0)[cc]{$b$}}
\put(69.,13.25){\makebox(0,0)[cc]{$c$}}
\put(13.75,16.75){\makebox(0,0)[cc]{$x_0$}}
\put(36,16){\makebox(0,0)[cc]{$x_0$}}
\put(57.75,15.5){\makebox(0,0)[cc]{$x_0$}}
\put(80.75,15.5){\makebox(0,0)[cc]{$x_0$}}
\put(50,3.75){\makebox(0,0)[cc]{$x_1$}}
\put(23.25,9.75){\vector(0,-1){9.5}}
\put(69,9.25){\vector(0,-1){9.5}}
\put(27,3.75){\makebox(0,0)[cc]{$x_2$}}
\put(72.25,3.75){\makebox(0,0)[cc]{$x_2$}}
\end{picture}
\vspace{1ex}

Picture~\theppp.
\end{center}

Notice that trees $T_1$, $T_3$ of the special forest have height $k$. So no carets can be placed over any pair of trees of the form $T_i$, $T_{i+1}$, where $0\le i\le3$. This explains why no edges labelled by $x_1^{-1}$, $x_2^{-1}$ can be accepted by vertices in the picture. The vertex $a$ correspons to a trivial marked tree so it does not accept $x_1$. However, it accepts $x_2$ since the tree to the right of $T_0$ has a caret. A similar situation holds for the vertex $c$. As for $b$, it accepts $x_1$ but it does not accept $x_2$ since the tree $T_2$ to the right of the marked tree $T_1$ is empty.

The condition on the tree $T_4$ allows us to avoid repetitions of vertices $a$, $b$, $c$ for different special forests.

Now we are subject to remove vertices of the form $a$, $b$, $c$ for all special forests. We also remove geometric edges incident to these vertices. For 3 vertices we thus remove no more than 14 directed edges from the graph. The key point here is inequality $\frac{14}3 < 5$. Here 5 is the limit of densities for $BB(n,k)$ considered as subgraphs in the left Cayley graph of $F$ in generators $\{x_0,x_1,x_2\}$. In the previous Section the same r\^ole was played by inequality $\frac{10}3 < 3{.}5$.

Now let us estimate densities of the new subgraphs $BB''(n,k)$ obtained from $B(n,k)$ after removing vertices and edges. We need to know the number of external edges labelled by $x_2$. The generating function for the set of vertices that do not accept $x_2$ is $\frac1{1-\Phi_k(x)}+\frac{x\Phi_k(x)}{(1-\Phi_k(x))^2}$. Here the marked tree is either rightmost, or the tree to the right of it is trivial. The number of such vertices is $\alpha_k(n)+\beta_k(n-1)$. The same number of vertices do not accept $x_2^{-1}$ due to Proposition~\ref{invlet} (symmetric property).

So we need to take all $6\beta_k(n)$ directed egdes incident to vertices of $BB(n,k)$ in the Cayley graph and subtract the number of external edges. This gives us $6\beta_k(n)-2\alpha_k(n)-2\beta_k(n-1)-2(\alpha_k(n)+\beta_k(n-1))=6\beta_k(n)-4\alpha_k(n)-4\beta_k(n-1)$ as the number of internal edges (see the number of external edges labelled by $x_0^{\pm1}$ and $x_1^{\pm1}$ in the previous Section).

The density of $BB(n,k)$ with respect to $\{x_0,x_1,x_2\}$ equals $6-4\frac{\alpha_k(n)}{\beta_k(n)}-4\frac{\beta_k(n-1)}{\beta_k(n)}\to6-4\xi_k$ as $n\to\infty$. Then it approaches $5$ as $k\to\infty$. 

To improve this estimate, we go to $BB''(n,k)$. We delete from $BB(n,k)$ no more than $14\sigma_k(n)$ directed edges and $3\sigma_k(n)$ vertices. The density of $BB''(n,k)$ has a lower bound 
$$
\frac{6\beta_k(n)-4\alpha_k(n)-4\beta_{k-1}(n)-14\sigma_k(n)}{\beta_k(n)-3\sigma_k(n)}=
$$
$$\frac{6-4\frac{\alpha_k(n)}{\beta_k(n)}-4\frac{\beta_k(n-1)}{\beta_k(n)}-14\frac{\sigma_k(n)}{\beta_k(n)}}{1-3\frac{\sigma_k(n)}{\beta_k(n)}}\to\frac{6-4\xi_k-14P(\xi_k)}{1-3P(\xi_k)}
$$
as $n\to\infty$. The latter approaches 
$$
\frac{5-14p}{1-3p}=5+\frac{p}{1-3p} > 5+\frac1{1200}
$$
as $k\to\infty$. So for some $n\gg k\gg1$ we can reach the density of finite subgraphs of the Cayley graph of $F$ in generators $\{x_0,x_1,x_2\}$ greater than $5{.}0008$. 
\vspace{1ex}

This completes the proof of Theorem~\ref{dens012}.
\vspace{1ex}

Notice that in~\cite{Gu21a} we thought of the generating set $A=\{x_0,x_1,x_2\}$ as a possible candidate to have {\em doubling property}. See~\cite[Section 8]{Gu22s} for general discussion of this property with respect to non-amenalility of groups. Briefly, this means that $|A^{\pm1}Y|\ge2|Y|$ for any finite nonempty subset $Y$ of the group $F$. This is equivalent to the fact that the outer boundary $\partial_oY$ has cardinality at least $|Y|$. Now we know this is not true. Moreover, even the Cheeger boundary $\partial_{\ast}Y$ may have cardinality strictly less than $|Y|$. This is exactly what Theorem~\ref{thm2} claims; it is a stronger fact because $|\partial_{\ast}Y|\ge|\partial_oY|$.

However, if we take a wider set of generators like $A=\{x_0,x_1,x_2,x_3\}$, then the probability of the automaton $BB(n,k)$ not to accept $x_i^{\pm1}$ is almost $\frac14$ for each $i\in\{1,2,3\}$. Even if we made a slight improvement of the density, then it is difficult to imagine that the isoperimetric constant will be far from $\frac32$. So we still have a chance that pure evacuation schemes on ${\cal C}(F;A)$ will exist for this case, according to~\cite{Gu22}. Constructing them will mean non-amenability of $F$.

\end{document}